\documentclass
[11pt,bezier]{article}
\usepackage{amsmath,amssymb,amsfonts,epsf,euscript}
\usepackage[usenames]{color}
\newpage

\newcommand{\be}{\begin{equation}}
\newcommand{\ee}{\end{equation}}

\renewcommand{\theequation}{\arabic{section}.\arabic{equation}}
\title{\bf\Large Estimation of Scale and Hurst Parameters of Semi-Selfsimilar Processes}
\vspace{-1cm}
\author
{S. Rezakhah\thanks{\scriptsize
Faculty of Mathematics and Computer Science, Amirkabir University of
Technology, 424 Hafez Avenue, Tehran 15914, Iran.
}\thanks{
School of Mathematics, Institute of Research in Fundamental Sciences (IPM),
PO Box: 19395-5746, Tehran, Iran.}
 \and  A. Philippe\thanks{
Laboratoire de Math\'ematiques Jean Leray, 2 rue de la houssinière, Universit\'e de Nantes, 44322 Nantes Cedex 3, France.
\,\,
E-mail:
namomath@aut.ac.ir(N. Modaresi), Anne.Philippe@univ-nantes.fr(A. Philippe), rezakhah@aut.ac.ir(S. Rezakhah).}
 \and  N. Modaresi$^*$}

\date{}
\begin{document}

\maketitle

\begin{abstract}
The characteristic feature of semi-selfsimilar process is the invariance of its finite dimensional distributions
by certain dilation for specific scaling factor. Estimating the scale parameter $\lambda$ and
the Hurst index of such processes is one of the fundamental problem in the literature.
We present some iterative method for estimation of the scale and Hurst parameters which is
addressed for semi-selfsimilar processes with stationary increments. This method is based on some
flexible sampling scheme and evaluating sample variance of increments in each scale intervals
$[\lambda^{n-1}, \lambda^n)$, $n\in \mathbb{ N}$.
For such iterative method we find the initial estimation for the
scale parameter by evaluating cumulative sum of moving sample variances and also by
evaluating sample variance of preceding and succeeding moving sample
variance of increments. We also present a new efficient method for estimation of Hurst parameter of selfsimilar processes. As an example we introduce simple fractional
Brownian motion (sfBm) which is semi-selfsimilar with stationary
increments. We present some simulations and numerical evaluation to
illustrate
the results and to estimate the scale for sfBm as a semi-selfsimilar process. We also present another simulation and show the efficiency of our method in estimation of Hurst parameter by
comparing its performance with
some previous methods.\\ \\
{\it Mathematics Subject Classification MSC 2010:} 62L12; 60G22; 60G18.
\end{abstract}

\section{Introduction}
Self-similarity has been discovered, analyzed and exploited in many frameworks, such as natural images \cite{b1-1}, fluctuations of stock market \cite{b2} and traffic modeling in broadband networks \cite{l1}.
The most frequently identified properties of high resolution traffic measurements {from a wide range of packet networks are long-range dependence and self-similarity \cite{m12}}. The fractal behavior of some important processes, such as fractional Brownian motion (fBm) and fractional Gaussian noise (fGn) is described by a single parameter, called Hurst parameter.

Self-similar processes are stochastic processes that are invariant in distribution under suitable scaling of time and space.
These processes { have increments which} enter naturally in the analysis of random phenomena (in time) exhibiting certain forms of long-range dependence \cite{e1}. Semi-selfsimilar or discrete scale invariant processes requires invariance by dilation for certain preferred scaling factors only.
In practice detecting self-similar property and estimation the self-similarity index and also scale parameter of the semi-selfsimilar processes are the main objects in the study of such processes. There are many problems arising in estimation of the Hurst parameter $H$ and determine those effects which can influence the results considerably. Estimation depends on several factors, e.g, the estimation technique, sample size, time scale, level shifts, correlation and data structure. Beran \cite {b2}, {Taqqu and Teverovsky \cite{taq1}, Bardet et.al. \cite{b1}} give a good review of  statistical aspects of parameter estimation methods for self-similar and long-memory processes. Among these the most well known are variance and covariance based methods (rescaled adjusted range R/S statistic and variogram including aggregated variance). The other methods are maximum likelihood based methods, Whittle estimator, the local Whittle \cite{r1}, { the local log-periodogram estimators, the global log-periodogram estimator and roughness-length method \cite{w2}.}
Traffic model based on the fBm contains three parameters: the mean rate, variance parameter and Hurst parameter.
Coeurjolly \cite{c2} developed a class of consistent estimators of the parameters of a fBm based on the asymptotic behavior of the $k$-th absolute moment of discrete variations of its sampled paths over a discrete grid of the interval $[0, 1]$. He derived explicit convergence rates for these types of estimators, valid through the whole range $0<H<1$. Consistent estimators of the fractal dimension of locally self-similar Gaussian processes based on convex combinations of sample quantiles of discrete variations of a sample path, almost sure convergence and the asymptotic normality for these estimators are derived \cite{c3}.

Several such studies of variations uncovered a generalization of fBm to non-Gaussian processes known as the Rosenblatt process and other Hermite processes. The processes are of order $q$ and with Hurst parameter $H\in (\frac{1}{2}, 1)$ which are self-similar with stationary increments and exhibit long-range dependence. The variations of these processes are studied and a consistent estimator for the self-similarity parameter from discrete observations of a Hermite process is constructed \cite{c1}.
For identifying a locally self-similar Gaussian process a new approach is proposed based on the asymptotic behavior of convex rearrangement, i.e., sums of increasing ordered variations \cite{p1}. They stated their result concerning the construction of the estimators of the local Holder index.
For the estimation of the Hurst parameter six methods are used that can be classified as temporal, spectral and time-scale methods, respectively.{ The temporal methods selected are (1) Rescaled range analysis (Beran \cite{b2}), (2) Level of zero crossings
(Coeurjolly \cite{c11}) and (3) Detrended fluctuation analysis (Hu et al. \cite{h1}). From the group of spectral methods, the (4) log-periodogram method (Beran \cite{b2}) is used while the (5) Wavelet transform modulus maxima method (Arneodo et al. \cite{a4})
and (6) Abry-Veitch estimator (Veitch and Abry \cite{v1}) are chosen from the group of the time-scale methods.}
Recently, moving average method to estimate the Hurst exponent and the correlation properties of the time series are presented and is compared with the rescaled range method of Hang Seng Index data for some periods \cite{w1}.

Of late there has been an increasing interest in describing and discussing the integer valued analogues of classical distributions, like stable, semi-stable, semi-selfdecomposable and geometrically infinitely divisible \cite{a3}, \cite{b4} and \cite{b5}.
The idea of binomial thinning to arrive at the right definitions of non-negative integer valued ($I_0$-valued) semi-selfsimilar processes are introduced \cite{s0}. In this class, they characterized semi-selfsimilar Levy processes in terms of an $I_0$-valued first order autoregressive series.

The paper is organized as follows. In section 2, we present flexible sampling method and some basic notions related to the self-similar and semi-selfsimilar processes in discrete parameter space. We introduce an example of semi-selfsimilar processes, named, simple fractional Brownian motion (sfBm) in this section too.
We present some simulations for sfBm with different scales and also simple Brownian motion with different Hurst indices as semi-selfsimilar processes in section 3.
In Section 4 we consider an iterative estimation method for scale parameter of semi-selfsimilar processes with stationary increments  based on sample variances of increments. For this iterative method we find the initial scale parameter by two procedures. In the first, we use the cumulative sum of moving sample variances (MSV) and in the second, sum of preceding and succeeding sample of variances of increments. We also present a new estimation method for estimating the Hurst parameter of selfsimilar processes based on  sample variance of different increments and show the efficiency of such new method by comparing its performance with some previous methods via simulation.

\renewcommand{\theequation}{\arabic{section}.\arabic{equation}}
\section{Method of flexible discrete sampling}
\setcounter{equation}{0}
We consider certain flexible sampling scheme as sampling at equally spaced points in each scale interval. Following this method of sampling
 from a semi-selfsimilar process $\{X(t), t\in \mathbb{ R^+}\}$ with scale $\lambda >1$, we decide to have some fix
  number of samples in each scale, say $T$.
  So we provide a discrete time semi-selfsimilar process $X(\cdot)$ with parameter space
 $\{ \lambda^{n-1} k\frac{(\lambda-1)}{T}, n \in \mathbb{ N}, k=0,1,\cdots T-1 \}$.\\

A process $\{X(t),t\in \mathbb{ R^+}\}$ is said to be self-similar of index $H>0$, if for every ${\lambda}>0$,
\be \{X(\lambda t)\}\stackrel{d}{=}\{\lambda^{H}X(t)\}\ee
where $\stackrel{d}{=}$ is the equality of all finite-dimensional distributions. As an intuition, self-similarity refers to an invariance with respect to any dilation factor. However, this may be a too strong requirement for capturing in situations that scaling properties are only observed for some preferred dilation factors.
This process is said to be semi-selfsimilar with index $H$ and scaling factor ${\lambda}_0>0$ if (2.1) holds just for $\lambda=\lambda_0^k$, {$k\in\mathbb{ N}$}.

A process $\{X(k),k\in {\check{T}}\}$ is called discrete time self-similar process with parameter space $\check{T}$,
where $\check{T}$ is any subset of countable distinct points of positive real numbers, if for any $k_1, k_2 \in \check{T}$
\vspace{-4mm}
\be \{X(k_2)\}\stackrel{d}{=}(\frac{k_2}{k_1})^H\{X(k_1)\}.\ee
The process $X(\cdot)$ is called discrete time semi-selfsimilar with scale $l>0$ and parameter space $\check{T}$, if for any $k_1, k_2=lk_1 \in \check{T}$, (2.2) holds; {see \cite{m11}}.
If the process $\{X(t),t\in \mathbb{ R^+}\}$ is semi-selfsimilar with scale $l>1$, then by the above method, one can find $\alpha$ by the equation $l=\alpha^T$, $T\in\mathbb{ N}$, then by sampling of the process at points $\alpha^{k}, k\in \mathbb{ Z}$, we have $X(\cdot)$ as a discrete time semi-selfsimilar process with parameter space $\check{T}=\{\alpha^{k}, k\in \mathbb{ Z}\}$.
Based on the definition of wide sense self-similar process presented in \cite{n1}, a random process $\{X(k),k\in \check{T}\}$ is called discrete time self-similar in the wide sense with index $H>0$ and parameter space $\check{T}$, where $\check{T}$ is any subset of distinct countable points of positive real numbers, if for all $k, k_1\in \check{T}$ and all $c>0$, where $ck, ck_1\in \check{T}:$\\

$(i)\,\,\ E[X^2(k)]<\infty$,

$(ii)\,\,E[X(ck)]=c^HE[X(k)]$,

$(iii)\,\, E[X(ck)X(ck_1)]=c^{2H}E[X(k)X(k_1)]$.\\ \\
If the above conditions hold for some fixed $c=c_0$, then the process is called discrete time semi-selfsimilar in the wide sense with scale $c_0$.
Through this paper we are dealt with the wide sense self-similar and wide sense scale invariant processes, where for simplicity we omit the term "in the wide sense" hereafter.\\
Using quasi Lamperti transform as ${\cal L}_{H,\alpha}Y(t)=t^HY(\log_\alpha t)$ \cite{m11}, the counterpart of any self-similar process $X(\cdot)$ can be identified as discrete time stationary process
$$Y(n)={\cal L}^{-1}_{H,\alpha}X(n)=\alpha^{-nH}X(\alpha^n).$$
It is clear by this relation that if $X(\cdot)$ is a discrete time semi-selfsimilar process with scale $l=\alpha^T$, $T\in \mathbb{ N}$ and parameter space $\check{T}=\{\alpha^k, k\in \mathbb{ Z}\}$, then $Y(\cdot)$ is a discrete time periodically correlated process with period $T$ and parameter space $\check{T}=\{n, n\in \mathbb{ Z}\}$.

\subsection{Simple fractional Brownian motion}
As an example of a semi-selfsimilar process, we introduce a process $X(t)$ with index $H>0$ and scale $\lambda>1$ which is called simple fractional Brownian motion (sfBm) and is defined as
\be\label{aa} X(t)=\sum_{n=1}^{\infty}\lambda^{(n-1)(H-H')}I_{[\lambda^{n-1}, \lambda^n)}(t)B_{H'}(t)\ee
where $B_{H'}(\cdot)$, $I(\cdot)$ are fractional Brownian motion with index $H'$ which is selfsimilar process and indicator function respectively. For $H'=\frac{1}{2}$, $B_{H'}(t)$ is Brownian motion and $X(t)$ is called simple Brownian motion (sBm) which has Markov property \cite{m11}.
It is easy to find that $X(t)$ is semi-selfsimilar with scale parameter $\lambda$ but not Markov.


\section{Simulation}

In this section we have generated and plotted simple Brownian motion (sBm) defined after equation (\ref{aa}) with scale $\lambda=1.2$ and Hurst indices $H=0.3,0.5,0.7$ at points
 $\lambda^n k(\lambda -1)/T$ for $T=20$ and $n=1,\ldots M$, and $k=0,1,\cdots T-1$where $M=20$.
Note that the fractional Brownian motion with index $H'$ is simulated using the circulant matrix embedding method (see \cite{b0} for a description).

We have samples of $M=20$ scale intervals
$[\lambda^n,\lambda^{n+1})$ and in each scale interval we have
$T=20$ equally space samples. Simple Brownian motion is a
semi-selfsimilar process with index $H$, and for $H=0.5$ it will be
Brownian motion, which has been plotted for
 comparison. By this method we show that how samples are effected by Hurst index.

We simulate and sfBm as a semi-selfsimilar process with $H=0.9$, and Hurst index
of the fBm as $H'=0.7$ for different scale parameter $\lambda=2, 4, 8$.

\input{epsf}
\epsfxsize=3in \epsfysize=2.2in
\begin{figure}
\centerline{$\hspace{-.1in}$\epsffile{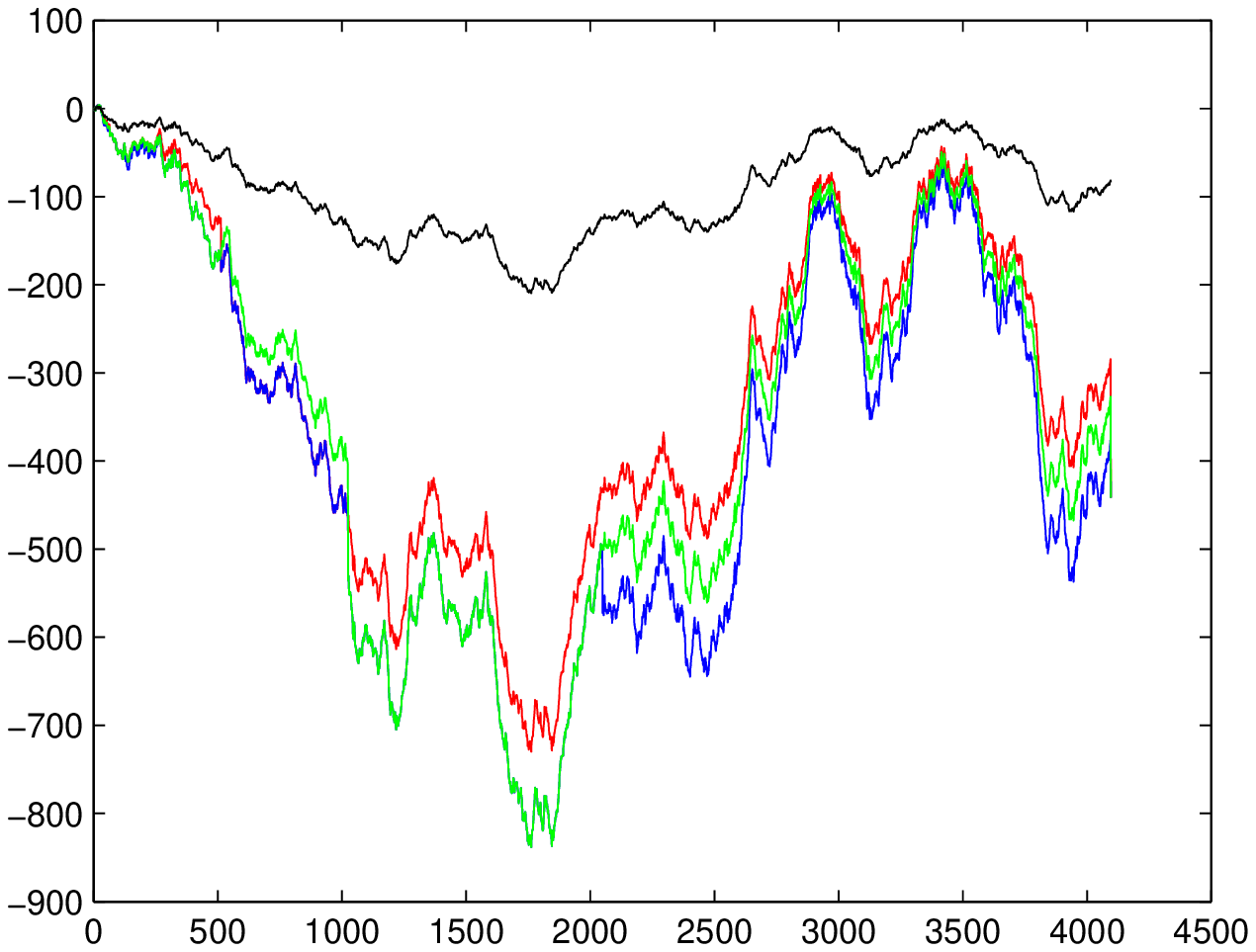}$\hspace{-.22in}$
\epsfxsize=3in \epsfysize=2.2in \epsffile{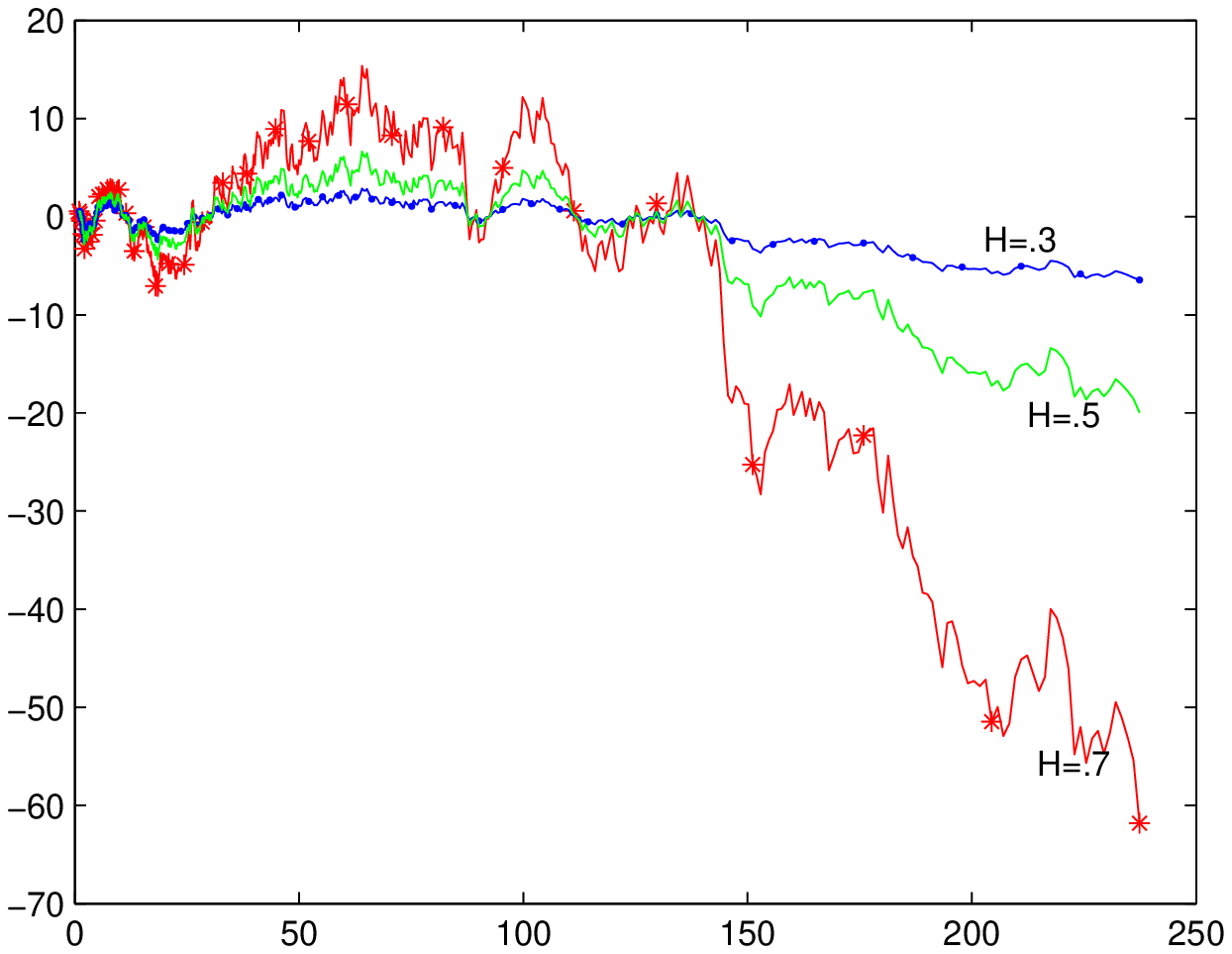}}
\caption{\scriptsize sfBm with Hurst indices H'=0.9 and H=0.7 and scales 2 [blue],4 [green] ,6 [red] and fBm with Hurst indice H=0.7 [black].  (Left),
sBm with scale 1.2 and Hurst indices 0.3 [blue], 0.5 [green], 0.7 [red] (Right), both with geometric
sampling and equally spaced in each scale interval}
\end{figure}

Figure 1 (Left) is the plotted sfBm which shows how the process
is effected by scale index. Both sBm and sfBm have semi-selfsimilar
property with stationary increments. In Figure 2, we show how our
method for recognizing scale parameter works. We have simulated sfBm
at equally space points with $H=0.9$ and $H'=0.2$ corresponding to
equation (2.1) for scale $\lambda=2$ and scale $\lambda=4$ left and
right correspondingly. We also evaluated and plotted corresponding
moving sample variance (MSV) of some small number of observations, say $b^*=10$, which is defined as
\be V_i=\frac{1}{b^*-1}\sum_{j=1}^{b^*}\big(Y(t_{i+j})-\bar{Y}(t_i)\big)^2,\ee
where $\bar{Y}(t_i)=\frac{1}{b^*}\sum_{j=1}^{b^*}Y(t_{i+j}), $ for $i=0, 1, \ldots, n-b^*$, and
 $n$ is the number of observations, in Figure 2. This figure shows how the sample variances can be
clustered for scale intervals and have jump at starting point of
each scale intervals. {Then the mean of the ratio of the length of successive
scale intervals} provides the initial estimate of scale parameter for
the iterative estimation method which is described in section 4.1.
\input{epsf}
\epsfxsize=3in \epsfysize=2.2in
\begin{figure}
\centerline{$\hspace{-.1in}$\epsffile{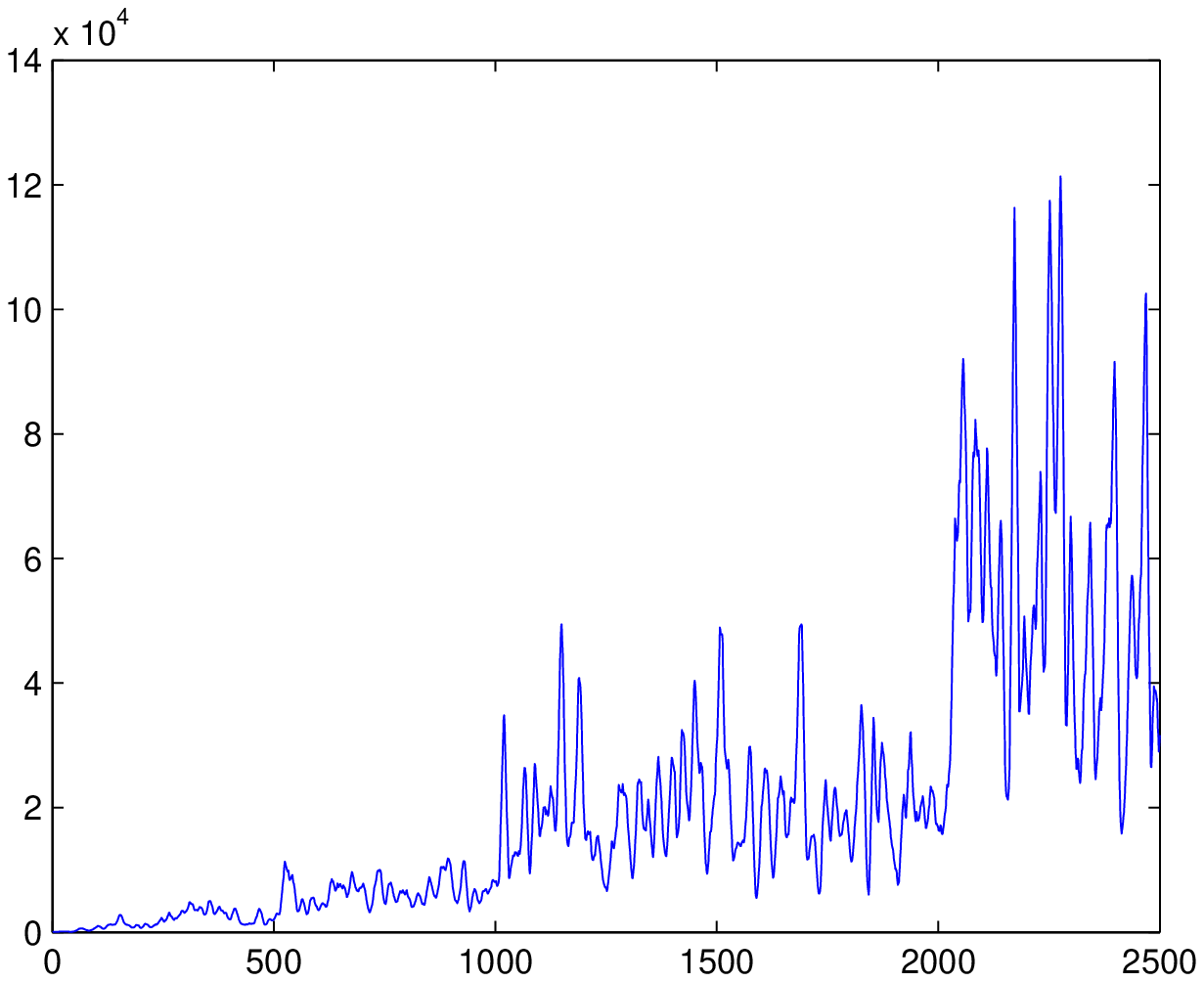}$\hspace{-.22in}$\epsfxsize=3in
\epsfysize=2.2in \epsffile{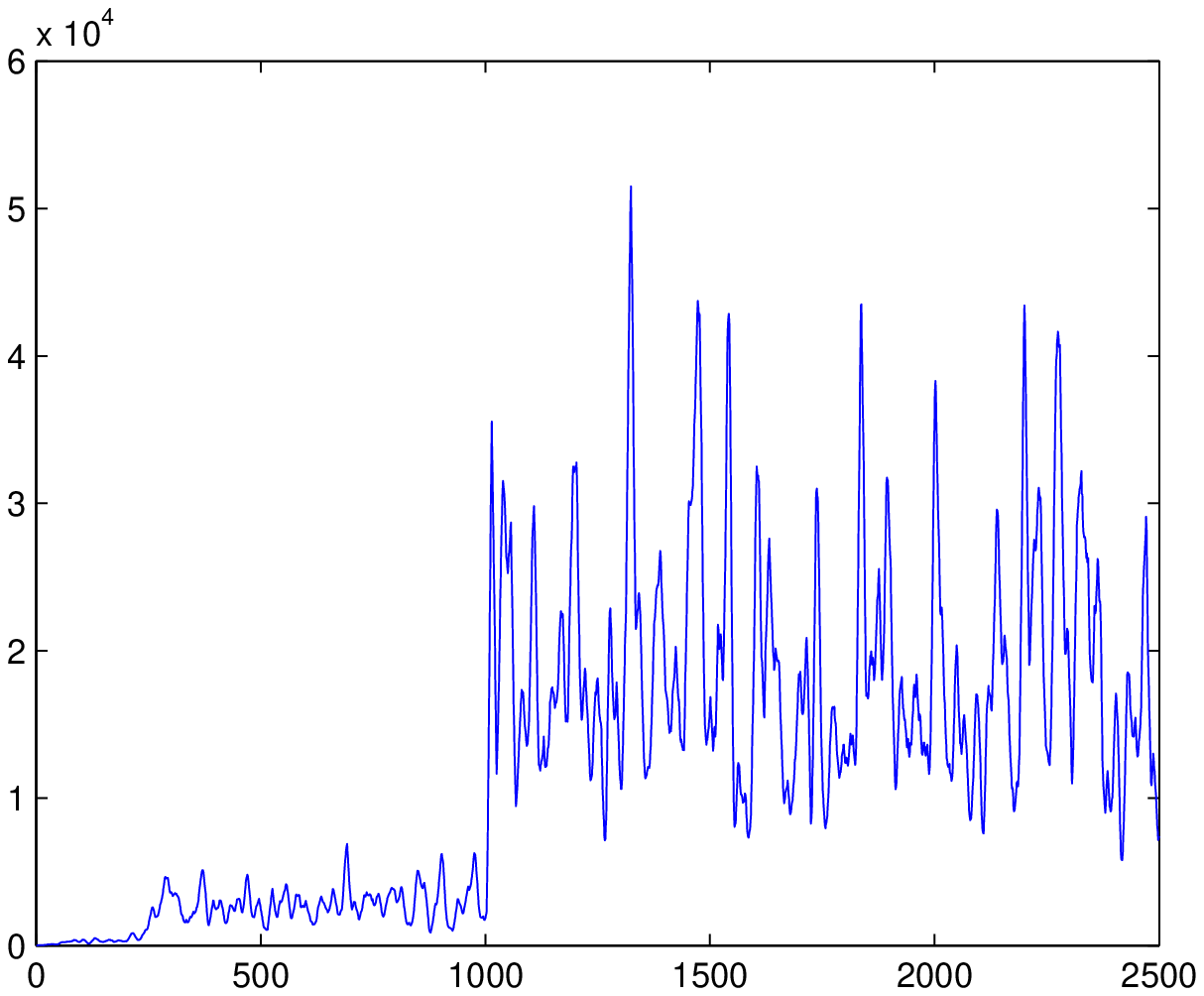}} \caption{\scriptsize
Moving average Sample variances with equally space sampling of sfBm,
Left figure Scale is 2 and right figure scale is 4.}
\end{figure}

\renewcommand{\theequation}{\arabic{section}.\arabic{equation}}
\section{Estimation procedure of parameters}
\setcounter{equation}{0}
As semi-selfsimilar processes are characterized by certain scale parameter $\lambda$ and Hurst index $H$. For the estimation of scale, one need to consider samples of several scale intervals $(\lambda^{k-1}, \lambda^{k}]$, $k\in\mathbb{ N}$.
For Estimation of the scale and Hurst index, we present a heuristic iteration method in this section which is proposed to estimate the scale parameter first and then the Hurst index. Our method can be applied in the situations where the process has stationary or independent increments. For such iteration method we need some starting points which has been described in section1.

\subsection{Initial choice of scale parameter}
{The starting point of our estimation method} can be evaluated by the following methods.
Suppose some equally space sampling of the process with small enough space to have a large number of samples, say $N$.
Using (3.4), one could evaluate moving sample variances $V_i$ for some small number, say $b^*$, of samples and plot $V_i$ against $i$, for $i=0,1, \ldots, N-b^*$, see  Figure 2.
This plot provide some clusters for the sample variance where jumps occurs at the starting point of all scale intervals $[\lambda^{k-1}, \lambda^{k})$, $k\in\mathbb{ N}$.
To cause these clusters to be distinguished well enough, we consider moving average of sample variance $V_i$'s  defined by (3.4) as
\be\label{W} W_i=\sum_{j=i}^{i+d^*}V_j/d^*,\ee
for some $d^*\in \mathbb{ N}$ say $d^*=20$, and $i=1,2,\cdots n^*$ where $n^*=n-b^*-d^*$, $b^*$ defined by (3.4).
By choosing appropriate values for $b^^*$ and $d^*$, and plotting consecutive samples of $W_i$'s one could
provide a clear distinguished clusters for samples of consecutive scale intervals.
By detecting start points of these consecutive scale intervals, one could estimate an appropriate initial value for the scale parameter $\lambda$.
   By calculating the cumulative sum  ( $W_i$'s ) defined by (\ref{W})
as $U_1, \ldots, U_N$, where $U_i=\sum_{j=1}^{i}W_j$,
as it is shown for example 1, by Figures 5(a) and 5(d), the plot of the sequence $U_i$ provide some broken lines that
broken points occurred at the border of the corresponding  clusters, or scale intervals at Figures 4(c) and 4(f).
So we can detect the starting points of scale intervals as these points and by these, the scale parameter
of a semi-selfsimilar process.
We consider following methods for such estimation.\\\

\noindent
{\bf First Method:}\\
 We can detect the starting points of successive scale intervals as change points on the mean of the sequence $W_i$.
Let $\tau_i, \; i=1,2,3$  to be the last three change points, or equivalently starting points of the last three scale intervals. So the initial choice of the scale parameter, say $\lambda_0$, can be evaluated as $\dfrac{\tau_{3}-\tau_{2} }{\tau_{2}-\tau_{1}}$.\\
Many methods for finding single or multiple change points  are available and  implemented in  the statistical softwares like {\bf R} with the package "changepoint".
Since the amplitude of the jump and the variance of the series increase, the methods for finding multiple change points
do not perform well. Therefore  we follow an iterative method for detecting change points in our special data.  At the first stage, we apply  a method for finding a single  change point on our series.  We use the non parametric method based on   the cumulative sums  test statistic (see \cite{chp} for details).   The method returns  the single most probable among  all possible changepoint locations.
Then  we delete the data from this point to the end of the data, and again we  iterate  the method in the remaining data.\\

\noindent
{\bf Second Method:}  \\
This method is based on calculation of sample variances of $W_i$s, defined by (\ref{W}),
 for $i=1,2\cdots z$ and for $i=z+1,\cdots n^*$ as sample variance
of those $W_i$s preceding some point $z$  as $L(z)$ and those succeeding point $z$ by $U(z)$,   as
$$L(z)=\frac{1}{z}\sum_{i=1}^{z}(W_i-\bar{W}_{1,z})^2,\hspace{1cm}U(z)=\frac{1}{n^*-z}\sum_{i=z+1}^{n^*}(W_i-\bar{W}_{2,z})^2,$$
$$\bar{W}_{1,z}=\frac{1}{z}\sum_{i=1}^{z}W_i,\hspace{1cm}\bar{W}_{2,z}=\frac{1}{n^*-z}\sum_{i=z+1}^{n^*}W_i,$$\\
 for all  $z=l^*,l^*+1, \ldots, n^*-l^*$ where $n^*$ defined by (\ref{W}), and for some appropriate $l^*\in\mathbb{ N}$, say $l^*=30$.
Then we evaluate
\begin{equation}\label{S} S(z)=L(z)+U(z)\hspace{.5in} \mbox{for all} \hspace{.5in} z=l^*, l^*+1, \ldots, n^*-l^*
\end{equation}   and plot them with respect to time $z$.
As it is shown in Figure 4, the minimum of this plot occurs at the starting point of the last scale interval, say $i_1$. By omitting samples from the point $i_1-j^*$ to the end, and by repeating this technique for samples $W_1,W_2,\ldots W_{i_1-j^*}$, for some appropriate $j^*\in\mathbb{ N}$, say $j^*=50$ one could find the minimum of $S(k), S(k+1),\ldots S(i_1-j^*-l^*)$ at the starting point of the last scale interval of the remaining samples, say $i_2$. Again by omitting samples from the point $i_2-j^*$ to the end and by repeating this method for the samples $W_1,W_2,\ldots W_{i_2-j^*}$,  the starting of the last scale intervals appears as the minima of the remaining samples, as $i_3$.
Then we evaluate the initial estimation of scale parameter, say $\lambda_0$ as $\frac{i_{1}-i_{2}}{i_{2}-i_{3}}$.

\subsection{Estimation of scale parameter}

Assume that $\{X(t), t \geq 1 \}$ is a discrete scale invariance (DSI), also called semi-selfsimilar,  process with some unknown scale $\lambda$. {  We present an iterative }  method for the estimation of the scale parameter $\lambda$ and  Hurst parameter $H$ for the case that the process has stationary increments.
Successive Scale intervals for such process is considered as $[\lambda^{k-1}, \lambda^{k})$, where $k\in \mathbb{ N}$.  We also assume that the process could have self similar property with some prescribed Hurst index $H'$ for samples inside each scale interval.
Proposing  equally space sampling implies that the increments have the same distribution. The iterativel estimation method is described for sfBm described in section 2.1  by the following steps.\\

\noindent
1- The whole duration of the study of the process is considered as the time interval [1,C].\\

\noindent
2- Sampling of the  process is assumed  at points $t_0, t_1, \ldots, t_N$, where
$$t_i=t_{i-1}+\frac{C-1}{N}$$
for $i=1, 2, \ldots, N$, and $t_0=1$. \\

\noindent
3-  {For implementing this method, we needs  some initial
 value for the scale parameter, say $\lambda_0>1$, which is  evaluated  by a
 numerical method,  described in section 4.1.}\\

\noindent
4- Successive Scale intervals are considered  as $[\lambda_0^{k-1}, \lambda_0^{k})$, $k=1,\ldots, M-1$, where $M$, the number of scale intervals, is evaluated as the largest $l\in\mathbb{ N}$ where $\lambda_0^l\leq C$.\\

\noindent
5- Increments of the process  are denoted by    $Y(t_i)=X(t_{i})-X(t_{i-1})$, for $i=1, \ldots, N$.  As  $E[B_{H'}(t_i)]=0$, so sample variance of the increments in $k$-th scale interval for this initial scale value $\lambda_0$, $k=0, 1, \ldots, M-1$ can be written as  , $S_{k}^2(\lambda_0)$  can be written as
$$S_{k}^2(\lambda_0)=\frac{1}{n_k} \!\!\!\!\! \sum_{\;\;\;\;i=N_{k-1}+1}^{N_k}\!\!\!\!\!\!(Y(t_{i})-Y(t_{i-1}))^2=\frac{\lambda^{2k(H-H')}}{n_k}\!\!\!\!\!\!
\sum_{\;\;\;\;i=1+N_{k-1}}^{N_k}\!\!\!\!\!\!\!\!
\big(B_{H'}(t_i)-B_{H'}(t_{i-1})\big)^2$$
where $N_k=\sum_{i=1}^{k}n_i$, $N_0=0$,  and $n_i$ is the number of samples in $i$-th scale interval.\\

\noindent
6- {For this iterative method}, several  other points, say more than 50 points, are considered in each side of $\lambda_0$ and at some small equally distances of each other. Assume that the total number of such points round initial estimate $\lambda_0$  is denoted by $m$. Rename these ordered points as $a_1,a_2,\cdots a_m$.
For all these points steps 4 and 5 are followed. Then for each $a_i$ we evaluate sum of the   sample variances corresponding to some of the  last scale intervals evaluated in step 5, say $j$ number, which covers at leat $\%95$ of observations  as
\begin{equation}\label{SS}
R(a_i)=\sum_{k=M-j}^M S^2_k(a_i).
\end{equation}
We plot such  $R(a_1),\ldots, R(a_m)$, which provide two clusters, those $R(a_i)$'s corresponding to the $a_i$'s less than true scale parameter, and those corresponding to $a_i$'s greater than true scale parameter.
In  example 1, such clusters are clearly shown by Figures 5(b) and 5(e).
Then for detecting change point of such identified  clusters and estimate true scale, we evaluate  sample variance of those $R(a_i)$ preceding each prescribed points $a_k$ as ${L}^*(a_k)$ and sample variance corresponding to those $R(\lambda_i)$ succeeding such prescribed $a_k$ as ${U}^*(a_k)$, for all $a_{k^*},\ldots, a_{m-k^*}$, for some appropriate $k^*$, say $k^*=20$. Thus
$$L^*(a_k)=\frac{1}{k-1}\sum_{i=1}^k (R(a_i)-\bar{R}_1)^2, \hspace{.2in}U^*(a_k)=\frac{1}{m-k-1}\sum_{i=k+1}^m (R(a_i)-\bar{R}_2)^2.$$
where $k=k^*,\ldots, m-{k^*}$, in which $\bar{R}_1$ and $\bar{R}_2$ are corresponding sample mean of $R(a_1),\ldots, R(a_k)$ and of
$R(a_{k+1}),\ldots, R(a_m)$ respectively.
Then we evaluate and plot
\begin{equation}\label{V}
{ V}(a_k )={ L}^*(a_k)+{ U}^*(a_k)
\end{equation}
for these points.
The point where $V(\cdot )$ has a minimum, which is considered as proper estimation of scale parameter is called $\lambda^*$.
In example 1, plots 5(c) and 5(f) show the minimum of such  $V(\cdot )$  for case 1 and case 2.\\

\noindent
7- Evaluate $\mu_k=\frac{S^2_{k}(\lambda^*)}{S^2_{k-1}(\lambda^*)}$, for $k= M-j,\ldots, M$.\\

\noindent
8- Use the equation $\bar{\mu}^*={(\lambda^*)}^{2(H-H')}$, to estimate $H-H'$, where $\bar{\mu}^*$ is the weighted average of  ${\mu}_k;$ $\; M-j\leq k\leq M-1$ evaluated at step 7, as
$$\bar{\mu}^*=\sum_{k=M-j}^M (\lambda^*)^{k-(M-j)}\mu_k/\sum_{k=M-j}^M (\lambda^*)^{k-(M-j)}.$$
{As one could estimate  $\lambda^*$  by step 6, this relation can be used  to estimate  ${H-H'}$}. So it remains to estimate $H'$, the Hurst index of underlying fBm. By dividing samples of the $k$-th scale interval to ${(\lambda^*)}^{k(H-H')}$ for $k=1,2,\ldots$  successive observations of fBm can be evaluated. Even though one can follow the usual methods for estimation the Hurst index of self-similar processes, we present a new method for such estimation by the followings.

\noindent
\subsection{ Estimation of Hurst parameter for H-sssi processes}
In this section we present two new methods for estimating Hurst parameter of self-similar process with stationary increment (H-sssi), which can be followed to estimate Hurst index $H'$ of the sfBm.
Let $\{X_{i}, i=1,2,\ldots, N\}$ be equally spaced samples of some H-sssi process as the main samples.
We consider some sub-samples at points $\{X_{i.k}, i=1,2,\ldots, [N/k]\}$  as the $k$-th sub-sample for
some fixed $k\in\mathbb{ N}$. The choice of the values of $k$ depends on the sample size we take for instance $k\in \{1,...,K_{\max} \}$.

For every $k\in \{1,...,K_{\max} \} $ we consider two sub-samples as $\{X_{i}\}$ and $\{X_{i.k}\}$, where  $i=1,2, \ldots, [N/k]$.
For these sub-samples, first and second order increments are defined as $Y_{1,i}=X_{i+1}-X_i$ and $Y_{2,i}=X_{i+2}-2X_{i+1}+X_i$ and as $Y_{1,i\cdot k}=X_{(i+1)\cdot k}-X_{i\cdot k}$ and $Y_{2,i\cdot k}=X_{(i+2)\cdot k}-2X_{(i+1)\cdot k}+X_{i\cdot k}$ respectively.
 By the followings we introduce two different methods for estimation  Hurst index, which are designed to evaluate
 $k^{2H}$ by the ratio of the sample variances of the first and second order
 increments of such sub-samples respectively. The first method is accurate for $H<0.75$ and the second method is more accurate for $H\geq 0.75$.\\

\noindent
{\bf First Method:}\\
For every $k\in \{1,...,K_{\max} \} $, we  evaluate sample variances of first order increments of the above mentioned sub-samples.  So we calculate  sample variance of the first order increments as  $Y_{1,1}, Y_{1,2}, \ldots Y_{1, [N/k]-1} $,  as $S^2_{1,k,1}$.  Also sample variance of the first order increments  $Y_{1,k}, Y_{1,2\cdot k}, \ldots, Y_{1, ([N/k]-1)k}$ as $S^2_{1,k,2}$, and evaluate
$\hat{H'}_k$ by the relation
\begin{equation}
\frac{S^2_{1,k,2}}{S^2_{1,k,1}}=k^{2\hat{H'}_k}\label{eq:1}
\end{equation}
where
$$S^2_{1,k,2}=\frac{1}{[\frac{N}{k}]-2}\sum_{i=1}^{[\frac{N}{k}]-1}(Y_{1,i\cdot k}-\bar{Y}_{1,k,2})^2\stackrel{d}{=}\frac{k^{2H'}}{[\frac{N}{k}]-2}\sum_{i=1}^{[\frac{N}{k}]-1}
(Y_{1,i}-\bar{Y}_{1,k,1})^2
=k^{2H'}S^2_{1,k,1},$$
and
$$\bar{Y}_{1,k,2}=\frac{1}{[\frac{N}{k}]-1}\sum_{i=1}^{[\frac{N}{k}]-1}Y_{1,i\cdot k}\stackrel{d}{=}\frac{k^{H'}}{[\frac{N}{k}]-1}\sum_{i=1}^{[\frac{N}{k}]-1}
Y_{1,i}=k^{H'}\bar{Y}_{1,k,1}$$
By (\ref{eq:1})  we find $\hat{H'}_k$ and finally estimate $H'$ as the mean of different $\hat{H'}_k$ which have been
 evaluated as:
$$\hat{H'}_k= \frac{1}{2(K_{\max}-1)}  \sum_{k=2}^{K_{\max}} \log\big(\frac{S^2_{1,k,2}}{S^2_{1,k,1}} \big) /\log(k).$$

\noindent
{\bf Second Method:}\\
We follow the same steps as mentioned in the first method, but this time based on the second order increments of sub-samples.
So we  evaluate sample variances of the second order increments  $Y_{2,1}, Y_{2,2}, \ldots Y_{2, [N/k]-2} $,  as $S^2_{2,k,1}$, and
 sample variance of the second order increments  $Y_{2,k}, Y_{2,2\cdot k}, \ldots, Y_{2, ([N/k]-2)k}$ as $S^2_{2,k,2}$, and evaluate
$\hat{H'}_k$ by the fact that
\begin{equation}
\frac{S^2_{2,k,2}}{S^2_{2,k,1}}=k^{2\hat{H'}_k}\label{eq:2}
\end{equation}
where
$$S^2_{2,k,2}=\frac{1}{[\frac{N}{k}]-3}\sum_{i=1}^{[\frac{N}{k}]-2}(Y_{2,i\cdot k}-\bar{Y}_{2,k,2})^2\stackrel{d}{=}\frac{k^{2H'}}{[\frac{N}{k}]-3}\sum_{i=1}^{[\frac{N}{k}]-2}
(Y_{2,i}-\bar{Y}_{2,k,1})^2
=k^{2H'}S^2_{2,k,1},$$
in which
$$\bar{Y}_{2,k,2}=\frac{1}{[\frac{N}{k}]-2}\sum_{i=1}^{[\frac{N}{k}]-2}Y_{2,i\cdot k}\stackrel{d}{=}\frac{k^{H'}}{[\frac{N}{k}]-2}\sum_{i=1}^{[\frac{N}{k}]-2}
Y_{2,i}=k^{H'}\bar{Y}_{2,k,1}$$
By (\ref{eq:2})  we find $\hat{H'}_k$ and finally estimate $H'$ as the mean of different $\hat{H'}_k$ which have been evaluated for different $k$ by this method as:
$$ \hat{H'}_k=\frac{1}{2(K_{\max}-1)}  \sum_{k=2}^{K_{\max}} \log\big(\frac{S^2_{2,k,2}}{S^2_{2,k,1}} \big) /\log(k).$$\\


\begin{figure}
\label{fig3}\vspace{0.1in}\vspace{.5in}

\centerline{$\hspace{-.2in}$\epsfxsize=3.2in \epsfysize=1.4in \epsffile{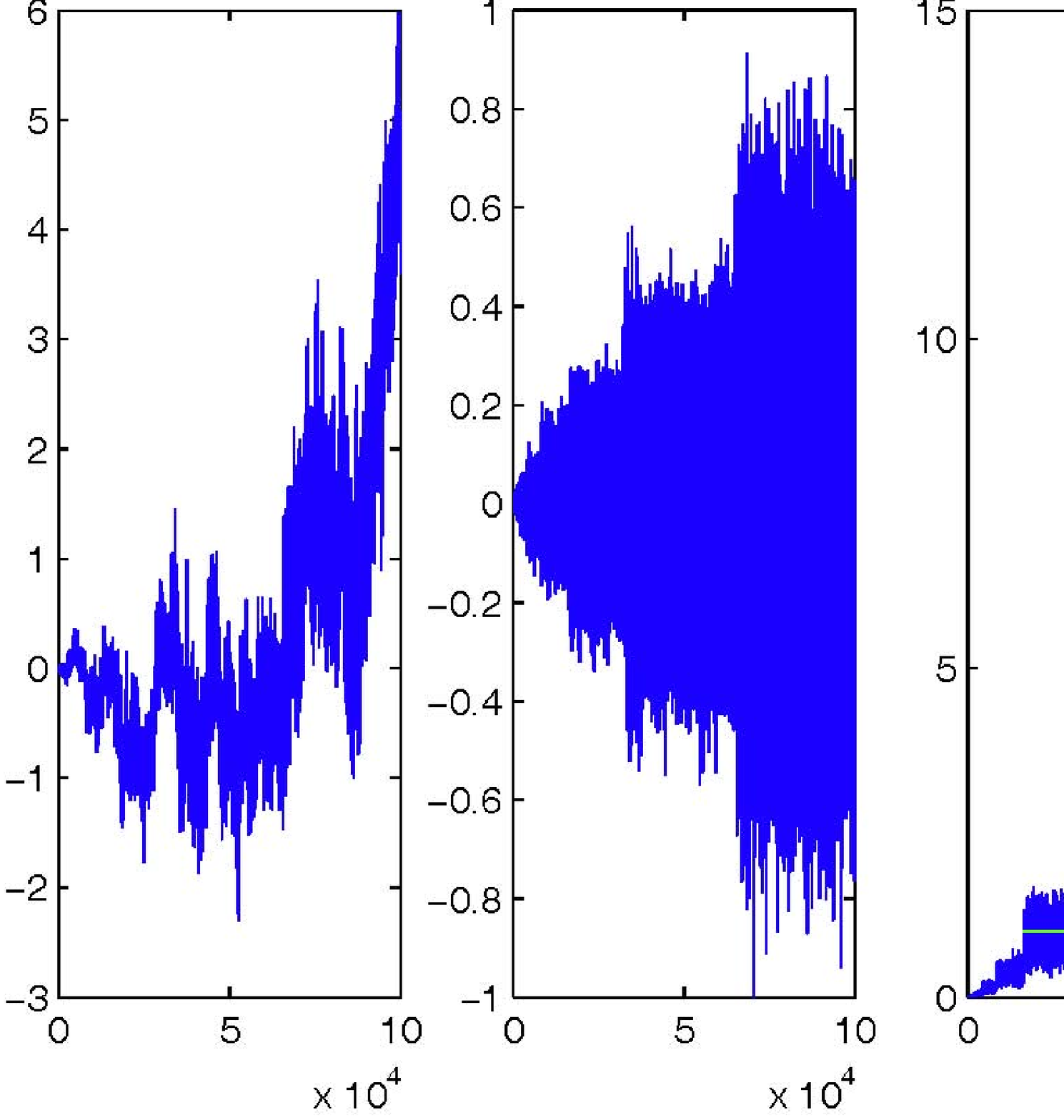}$\hspace{-.3in}$
\epsfxsize=3.2in
\epsfysize=1.4in
 \epsffile{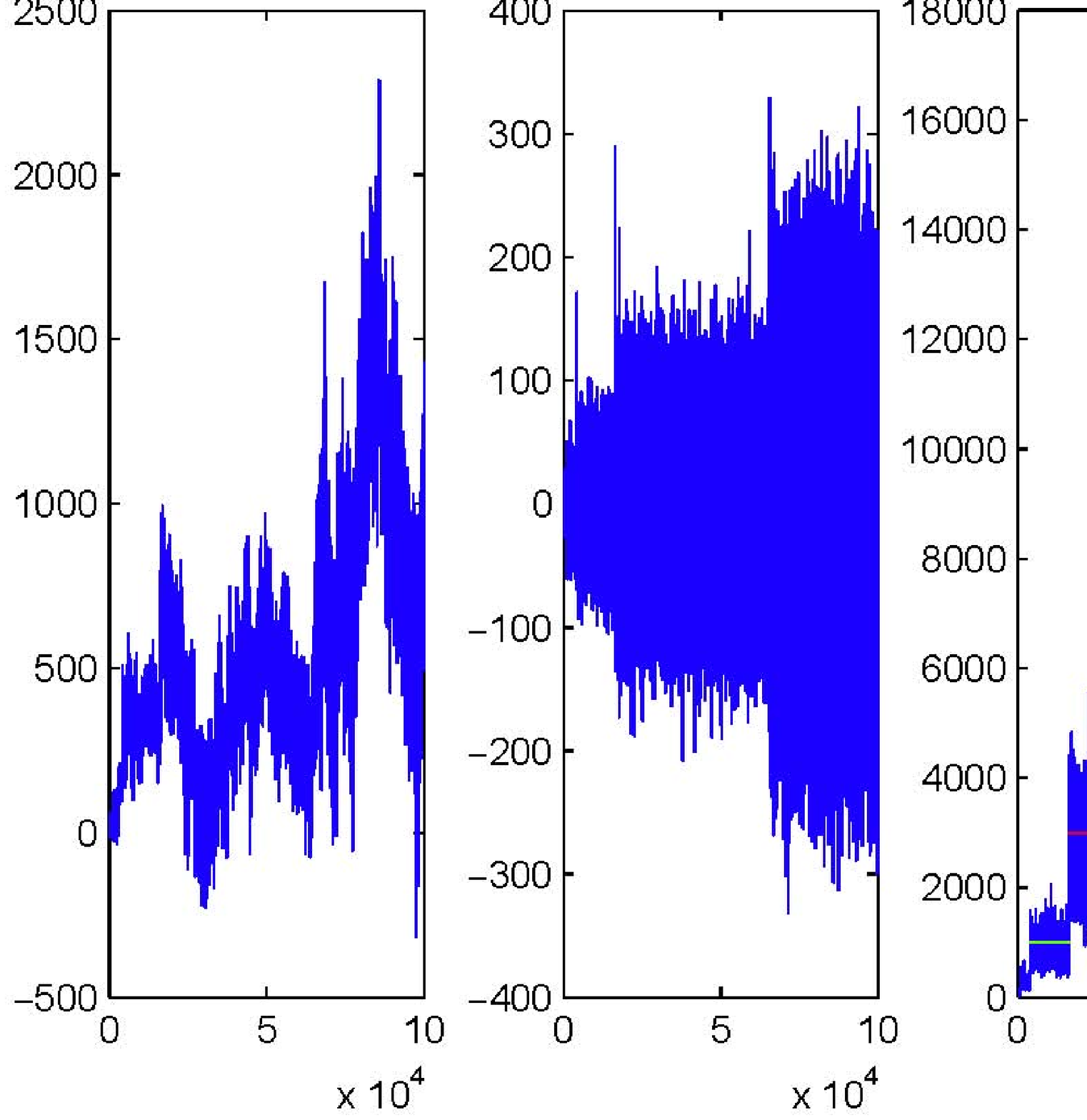}}\vspace{-.7in}

\caption{\scriptsize
Left three plots are 1-sfBm, 2-Increments of sfBm, and 3-Moving average $W_i$ of MSV of the increments of  sfBm, all with scale 2. The right three figure are corresponding figures with scale 4..}\vspace{0in}
\end{figure}
\vspace{0.2in}

\noindent
{\bf\large Example 1}\\
In this example we simulate sfBm processes defined by (2.3) with  scales $\lambda$ and Hurst parameters $(H,H')$ as semi-selfsimilar processes with stationary increments, and then apply our methods to estimate their parameters. For doing this we consider two cases first, $\lambda_1=2$, $H_1=0.9$, ${H'}_1=0.2$
and second, $\lambda_2=4$, $H_2=0.6$, ${H'}_2=0.2$ and simulate 100,000 samples in each case, which have been plotted in Figures 3 or plan is to estimate scales first. So we follow to estimate some initial value by method one and two, described in section 4.2.
Using the first method we  detect change points in the corresponding MSV plotted in Figure 3 at the starting points of last three scale intervals as $a_1=65480, b_1= 32745 , c_1=16354$ for the first case and as $a_2=65497, b_2=16349, c_2= 4060 $
for the second case. So the corresponding initial value of the scale parameter would be $\hat{\lambda_1}=1.99$ and $\hat{\lambda_2}=3.99$ by the first method.
%
%
%

\input{epsf}
\begin{figure}\vspace{.2in}

\epsfxsize=3.2in \epsfysize=1.7in
\centerline{$\hspace{-.1in}$ \epsfxsize=3.2in \epsfysize=1.7in
\epsffile{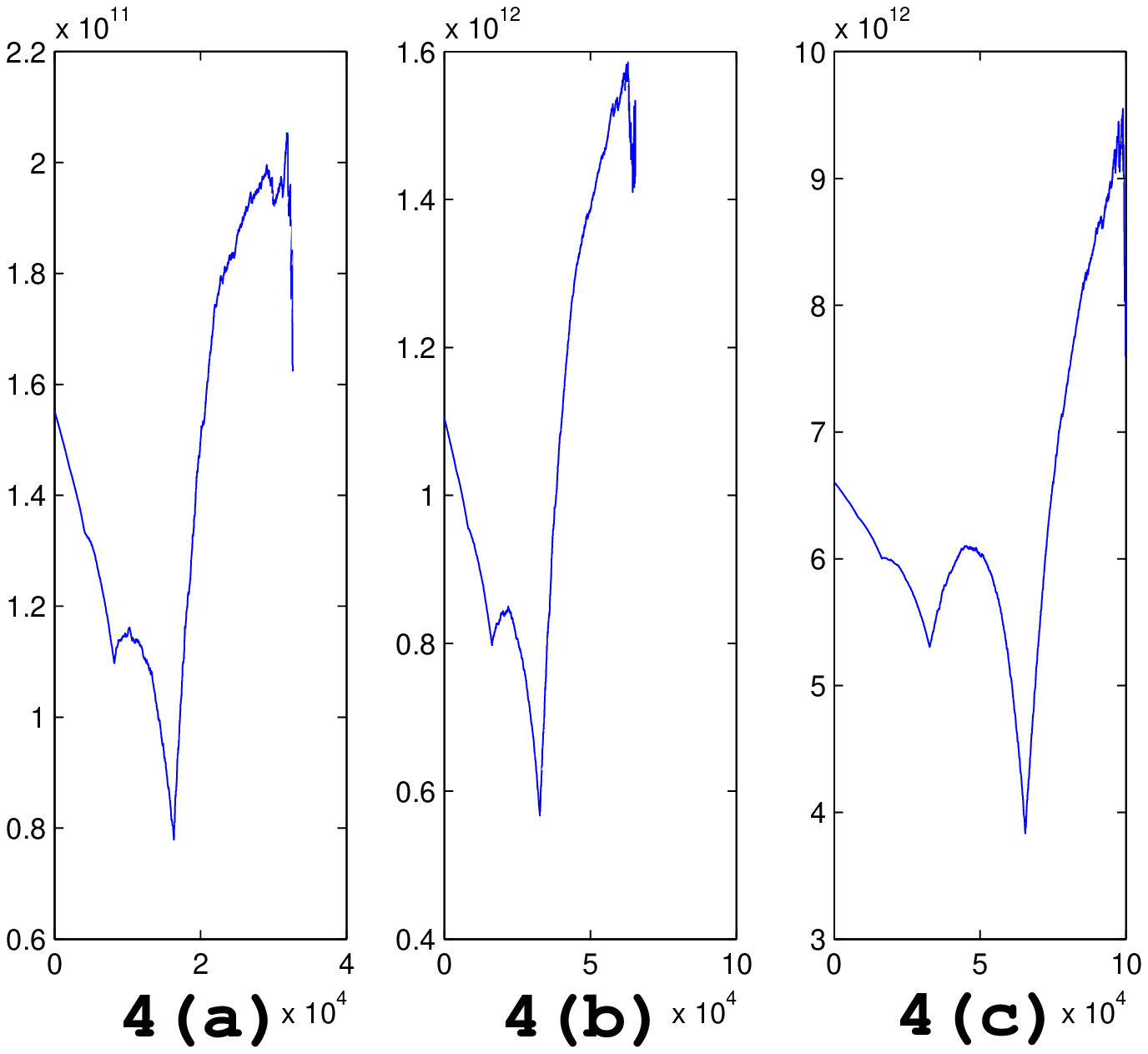}$\hspace{-.25in}$\epsfxsize=3.2in
\epsfysize=1.7in \epsffile{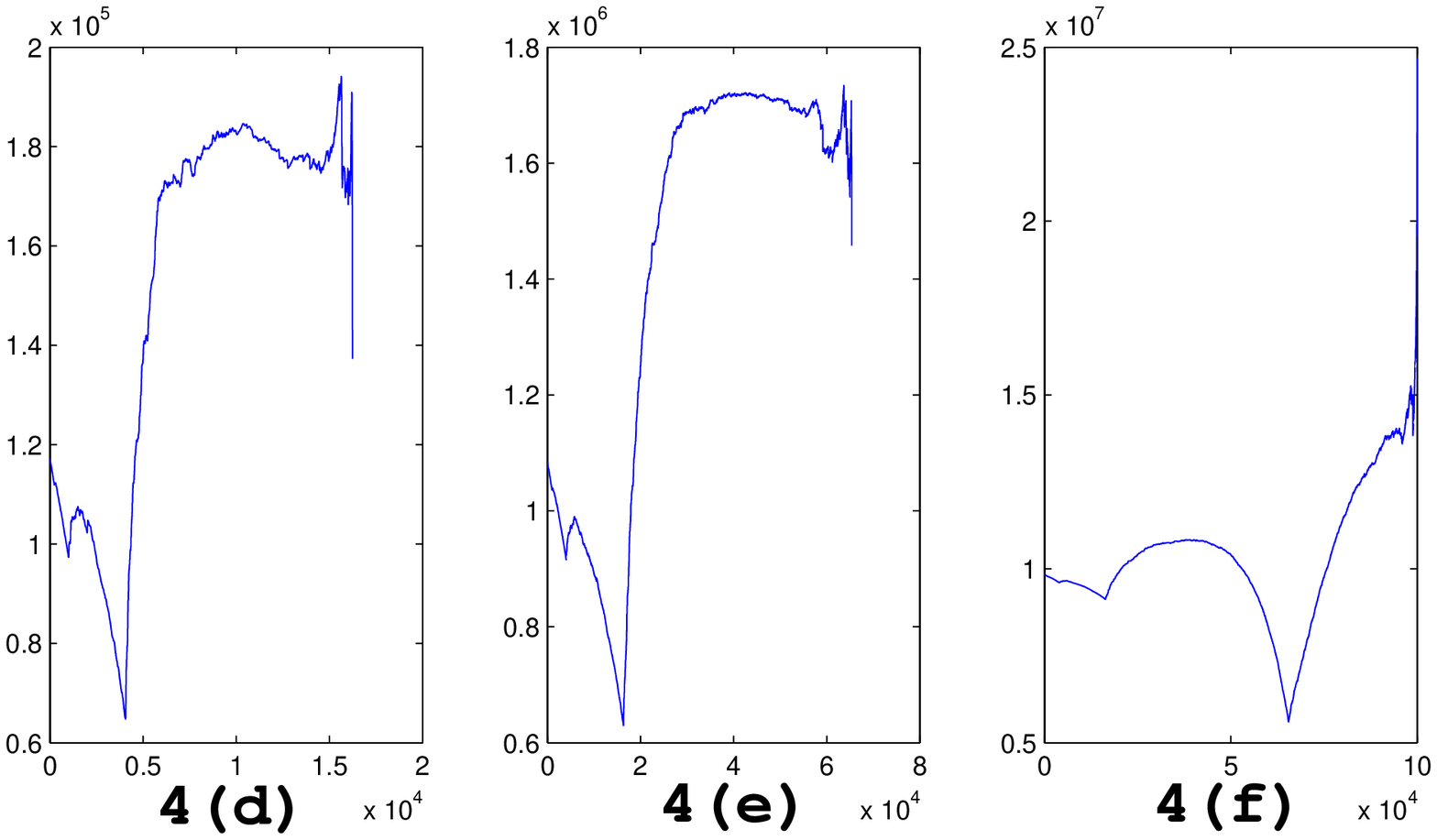}}$\vspace{.2in}$

\caption{\scriptsize
4(c)-Plot of Sum of sample variances $S(z)$, defined by (\ref{S}) for the whole samples $W_i,\; i=1,\cdots, n^*$, for case 1, example 1,which has minima at the start of last scale interval, say $i_1$.
4(b)- Plot of $S(z)$ but for remaining samples $W_i$, when samples from $i_3-l^*$ to the end of data is removed, which has minima at $i_2$, starting point of the last scale interval of remaining samples. 4(a)- Plot of $S(z)$, this time by removing samples from $i_2-j^*$ to the end of data, which has minima at $i_1$, starting point of the last scale interval of remaining samples. The last three figures are corresponding figures for case 2 of example 1.}
\end{figure}

%
%
%
%

\input{epsf}
\epsfxsize=3.2in \epsfysize=1.7in
\begin{figure}\vspace{.4in}

\centerline{$\hspace{-.1in}$
\epsffile{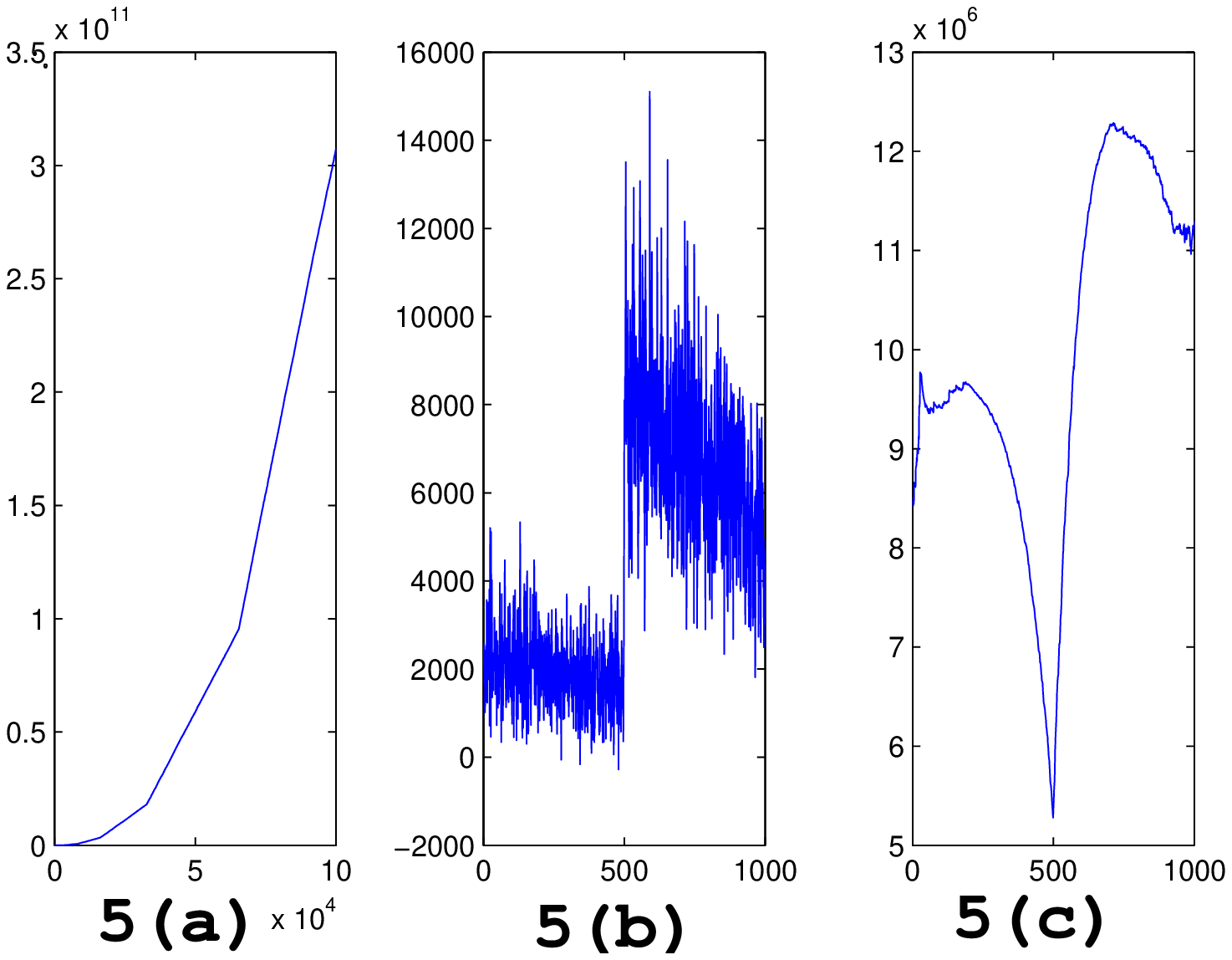}$\hspace{-.25in}$\epsfxsize=3.2in
\epsfysize=1.7in \epsffile{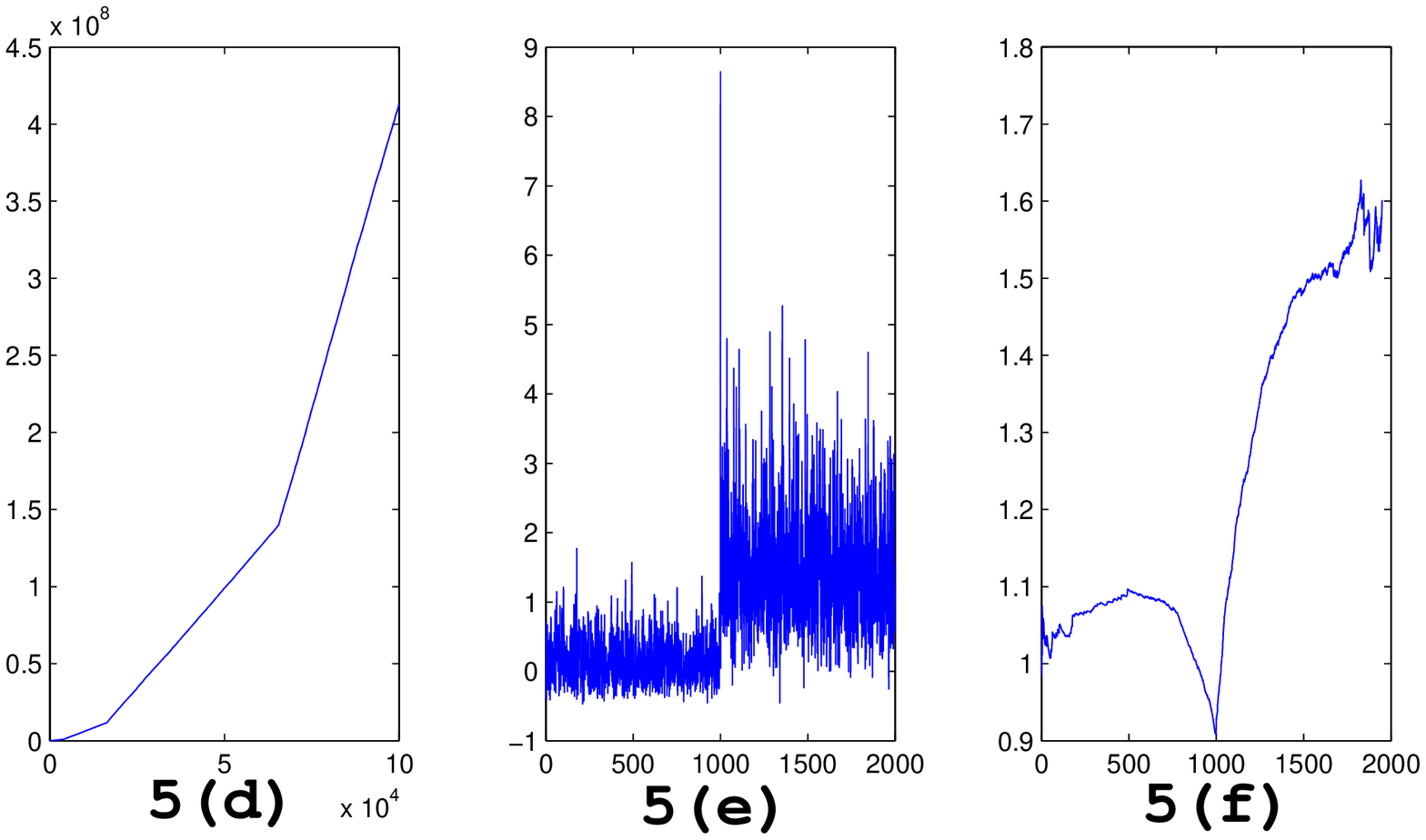}} \vspace{0.2in}

\caption{\scriptsize
5(a)-Cumulative sum of moving sample variances of increments of sfBm of case 1, example1, 5(b)-Sum of sample variances of increments of $R(a_i)$, defined by $(\ref{SS})$,  for $1000$ equally spaced samples in the interval $[1.9501,\cdots 2.05]$ round the initial estimate of scale $\lambda_0=1.9995$, and 5(c)- Sum of sample variances $V(a_k)$, defined by (\ref{V}) of preceding and succeeding samples to each point $a_k$ in Figure 4(b), which has minima at the change point of figure 4(b), $\lambda^*=2$, as the best estimate of scale parameter,  for case 1,example1. Right three figures are corresponding figures for case 2, example 1, that shows detection of the best estimator for scale as $\lambda^*=4$ in Figure 5(f).}
\end{figure} $\vspace{.5in}$

By the second method we detect the starting points of the last three scale intervals for case one as $65519,32743,16351$, which is shown by Figures 4(c), 4(b) and 4(a), which corresponding statistics $S(z)$ is defined by relation (\ref{S}). Thus initial values of scale parameter for case one  is evaluated as $\tilde{\lambda_1}=1.9995$.
For case 2 the starting points of the last three scale interval is estimated as  $65498,16345,4056$ which has been shown by Figure 4(d), 4(e) and 4(f), as the minima of $S(z)$ defined by (\ref{S}). Thus  initial value of scale parameter is evaluated as $\tilde{\lambda_2}=3.9998$. Then following the iterative iterative estimation method
we  evaluate $V(a_k)$ which defined by (\ref{V})  for $1000$ equally spaced points in the interval $[1.95, 2.05]$, round initial estimate $\lambda_0=1.9995$,  which has been plotted by 5(c), and  has minima at $\lambda^*_1=2$ as the best estimate for scale parameter in case 1,.  Also be evaluating $V(a_k)$ for $2000$ equally space points in the interval $[3.9, 4.1]$ round initial estimate $\lambda_0=3.9998$, which has been plotted by (5.f) and has minima at $\lambda^*_2=4$ as the best estimate for scale parameter in case 2.
We evaluated sum of sample variances  of five last scale intervals for case one and for three last scale intervals for case 2 with estimated scale parameters $\lambda_1^*=2$ and $\lambda_2^*=4$ respectively, so that our scale intervals cover at least $\%95$ of samples. Then we evaluated  $\bar{\mu}^*_1=2.6480$ and $\bar{\mu}^*_2=3.0108$ corresponding to $\lambda^*_1$ and  $\lambda^*_2$ respectively as described in step 8 of sequential estimation in section 4.1. Therefore we evaluate our estimations as $\hat{H}_1-\hat{H'}_1=0.7024$ and $\hat{H}_2-\hat{H'}_2=0.3992$.

Finally we extract the corresponding samples of underlying fBm as H-sssi process, by dividing samples in $k$-th
scale interval of sfBm in case 1  by ${(\lambda^*_1)}^{H_1-{H'}_1}= \bar{\mu}^*_1=1.6273$, so  ${H_1-{H'}_1}= 0.7025 $, and in case 2 by
$ {(\lambda^*_2)}^{H_2-{H'}_2}=\bar{\mu}^*_2=1.7352$, so ${H_2-{H'}_2}=0.3976$.  By applying our method described above in section 4.3, we find for case 1 that $ \hat{H'}_1=0.1957,\;\; \hat{H}_1=0.8981$  and for case 2. So we evaluate $ \hat{H'}_2=0.2039$ and $\hat{H}_2=0.6031$.\\\

\noindent
{\bf Efficiency of the Hurst parameter Estimation}\\\

\noindent
To visualize efficiency of  our Hurst estimation methods, we have simulated 10000 samples of fractional Brownian Motion with different Hurst parameter with 500 repetition  and  presented the graph of the mean square errors (MSE) for Hurst Estimation of our methods and also of the method of Quadratic Variation \cite{ist} and the method of Convex Rearrangement \cite{p1}. We have used the difference of order two for the last two methods.
As it is shown by Figure 6 for $H<0.75$ our first method has less mean square error in compare to other methods and so is the best. For $H\geq 0.75$ our second method
has less mean square error in compare with other methods and so is the best.
\input{epsf}
\epsfxsize=3in \epsfysize=2in
\begin{figure}
\centerline{$\hspace{-.1in}$\epsffile{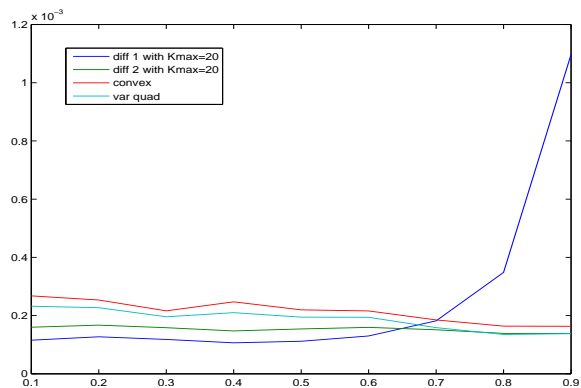}
}
\caption{\scriptsize
 Mean square error in Estimation of Hurst index of using 10000 samples of fBm with 500 repetition.
}
\end{figure}

%

\section{Conclusion}
 This paper provides some innovative method for estimation of  the scale parameter of some class of semi-selfsimilar processes. For this we consider some starting point for iterative estimation of scale parameter and then by some precise method we estimate it.
There are various open issues and vivid discussions about the estimation of scale parameter but still there is not a universal method which could be considered as the most promising method to find the best approximation of scale parameter in all cases.
 For the processes with stationary increments, considering equally space sampling implies that the increments have the same distribution inside each scale interval, so this enables us to  consider such an iterative estimation method based on moving sample variances.
Simulations and numerical evaluations clarified our results for the simple fractional Brownian motion.
By detecting start points of consecutive scale intervals, we obtained an initial value for the scale parameter.
Moving sample variances provide some clusters where jumps occurs at the starting point of scale intervals. Our iterative method provided a very good estimation for scale parameter.
 We also presented a heuristic method for the estimation of Hurst parameter of self-similar process which has less MSE in compare to previous methods.  Moreover our estimation methods are easily implemented and are computationally fast.

\bibliographystyle{plain}

\begin{thebibliography}{}

\bibitem{a2} P. Abry, P. Flandrin, M. Taqqu and D. Veitch (2000). Wavelets for the analysis, estimation and synthesis of scaling data, {\em Self-Similar Network Traffic and Performance Evaluation}, K. Park and W. Willinger, Eds. New York: Wiley.
\bibitem{a3} E. E. Aly and N. Bouzar (2000). On geometric infinite divisibility and stability, {\em Ann. Inst. Statist. Math.}, 52, 790-799.
\bibitem{a4} A. Arneodo, E. Bacry, P.V. Graves and J.F. Muzy (1995). Characterizing long-range correlations in DNA sequences from wavelet analysis, {\em Phys. Rev. Lett.}, 74(16), 3293-3296.
\bibitem{b0} J. M. Bardet, G. Lang, G. Oppenheim, A. Philippe and M.S. Taqqu (2003). Generators of long-range dependent processes: A survey. {\em Long-Range Dependence: Theory and Applications} (Eds P. Doukhan, M.S. Taqqu, G. Oppenheim),  Birkhauser, Boston, pp. 579-623.
\bibitem{b1} J. M. Bardet, G. Lang, G. Oppenheim, A. Philippe, S. Stoev and M.S. Taqqu (2003). Semi-parametric estimation of the long-range dependent processes: A survey. {\em Long-Range Dependence: Theory and Applications} (Eds P. Doukhan, G. Oppenheim and M.S. Taqqu), Birkhauser, Boston, pp. 557-577.
\bibitem{b1-1} M. Barnsley (1988). Fractals everywhere, {\em Academic Press, San Diego}.
\bibitem{b2} J. Beran (1994). Satistics for long memory processes, {\em Chapman and Hall, New York}.
\bibitem{b3} P. Borgnat, P. Flandrin and P.O. Amblard (2002). Stochastic discrete scale invariance, {\em IEEE SignalProcess}, Lett.9, pp.181-184.
\bibitem{b4} N. Bouzar (2008). Semi-self-decomposable distributions on $\mathbb{ Z^+}$, {\em Annals of the Institute of Statistical Mathematics}, Volume 60, Number 4, 901-917.
\bibitem{b5} N. Bouzar and K. Jayakumar (2008). Time series with discrete semi-stable marginals, {\em Statistical Papers}, Volume 49, Number 4, 619-635.
\bibitem{c1} A. Chronopoulouc, C. Tudor and F.G. Viens (2010).  Self-similarity parameter estimation and reproduction property for non-Gaussian Hermite processes, {\em Communications in Stochastic Analysis}.
\bibitem{c11} J.F. Coeurjolly (2000). Simulation and identification of the Fractional brownian motion: a bibliographical and comparative study, {\em J. Stat. Soft.}, 5(7).
\bibitem{c2} J.F. Coeurjolly (2001). Estimating the parameters of a fractional Brownian motion by discrete variations of its sample paths, {\em Statistical Inference for Stochastic Processes}, 4, 199-227.
\bibitem{c3} J.F. Coeurjolly (2008). Hurst exponent estimation of locally self-similar Gaussian processes using sample quantiles, {\em  Annals of Statistics}, 36(3), 1404-1434.
\bibitem{chp} M. Csorgo and L. Horvath (1997). \textit{Limit Theorems in Change-Point Analysis}, Wiley
\bibitem{e1} P. Embrechts and M. Maejima (2000). An Introduction to the Theory of Self-Similar Stochastic Processes, {\em International Journal of Modern Physics B}, Volume 14, Issue 12-13, pp. 1399-1420.
\bibitem{h1} K. Hu, P. Ivanov, Z. Chen, P. Carpena and H.E. Stanley (2001). Effect of trends on detrended fuctuation analysis, {\em Phys. Rev.}, E 64, 011114.
\bibitem{ist} J. Istas and G. Lang (1997). Quadratic variations and estimation of the local
    Holder index of a Gaussian process, {\em Ann. Inst. H. Poincar�e, Probab. Statist.}, Vol 33(4), pp.407-436.
\bibitem{l1} W.E. Leland, M.S. Taqqu, W. Willinger and D.V. Wilson (1994). On the self-similar nature of Ethernet traffic (extended version), {\em IEEE/ACM Transactions on Networking}, Vol.2, pp.1-15.
\bibitem{m1} D. Markovic and M. Koch (2005). Sensitivity of Hurst parameter estimation to periodic signals in time series and filtering approaches, {\em Geophysical Research Letters}, Vol.32, L17401 (5pp).
\bibitem{m11} N. Modarresi and S. Rezakhah (2010). Spectral analysis of Multi-dimensional selfsimilar Markov processes, {\em Journal of Physics A: Mathematical and Theoretical}, Vol.43, No.12, 125004 (14pp).
\bibitem{m12} S. Molnar, A. Vidacs and A. A. Nilsson (1997). Bottlenecks on the way towards fractal characterization of
network traffic: Estimation and interpretation of the Hurst parameter, {\em  International
Conference of the Performance and Management of Complex Communication Networks}, Tsukuba, Japan, PMCCN'97, pp.111-134.
\bibitem{n1} C. J. Nuzman and H. V. Poor (2000). Linear estimation of self-similar processes via Lamperti transformation, {\em Journal of Applied Probability}, No.37(2).
\bibitem{p1} A. Philippe and E. Thilly (2002). Identification of a locally self-similar Gaussian process by using convex rearrangement, {\em Methodology and Computing in Applied Probability}, No.4, 195-209.
\bibitem{r1} P.M. Robinson (1995). Gaussian semiparametric estimation of long range dep endence, {\em The Annals of Statistics}, 23, 1630-1661.
\bibitem{s0} S. Satheesh and E. Sandhya (2008). Semi-Selfsimilar Processes on $R$ and $I_0$, {\em International Mathematical Forum}, 3, No.36, 1781-1784.
\bibitem{taq1} M.S. Taqqu and V. Teverovsky (1998). On estimating long-range dependence in finite and infinite variance series, {\em A Practical Guide to Heavy Tails: Statistical Techniques and Applications}, R.J. Adler, R.E. Feldman and M.S. Taqqu, editors, Birkhauser, pp. 177-217.
\bibitem{v1} D. Veitch and P. Abry (1999). A wavelet-based joint estimator of the parameters of long-range dependence, {\em IEEE Trans. Inf. Theory 45}, 878897.
\bibitem{w1} N. Wang, Y. Li and H. Zhang (2010). Hurst exponent estimation based on moving average method, {\em Advances in Wireless Networks and Information Systems}, Vol.72, 137-142.
\bibitem{w2} J. Wawszczak (2005). Methods for estimating the Hurst exponent. The analysis of its value for fracture surface research, {\em Materials Science-Poland}, Vol.23, No.2.
\end{thebibliography}

\end{document}